\documentclass[10pt]{article}
\usepackage{amsmath}
\usepackage{amsfonts}
\usepackage{latexsym}
\usepackage{amssymb}
\usepackage{enumerate}
\usepackage[dvips]{graphics}

\newcommand{\Z}{{\mathbf Z}}
\newcommand{\Q}{{\mathbf Q}}
\newcommand{\R}{{\mathbf R}}
\newcommand{\C}{{\mathbf C}}

\newcommand{\Hom}{{\rm Hom}}

\newcommand{\qed}{\begin{flushright} $\Box$\end{flushright}}
\newcommand{\supp}{\rm {supp}}

\newcommand{\comment}[1]{}

\def\U{\mathcal U}

\def\V{\mathcal V}

\def\ker{{\rm Ker }}

\def\lra{\longrightarrow}

\def\bbr{{\mathbf R}}
\def\bbz{{\mathbf Z}}

\def\cu{{\mathcal U}}
\def\cv{{\mathcal V}}
\def\x{\times}

\def\eps{\varepsilon}
\def\del{\partial}
\def\pHi{\varphi}

\def\R{\mathbf R}

\newtheorem{theorem}{Theorem}

\newtheorem{prop}{Proposition}
\newtheorem{lemma}{Lemma}

\newtheorem{corollary}{Corollary}
\newtheorem{definition}{Definition}

\numberwithin{equation}{section}
\newenvironment{proof}{\par\noindent{\bf Proof.} }{\par\noindent}

\title{Lyapunov 1-forms for flows}
\author{M.~Farber\thanks{Partially supported by a grant from the Israel Academy of Sciences and Humanities
and by the Herman Minkowski Center for Geometry; part of this work
was done while M. Farber visited FIM ETH in Zurich},\,
T.~Kappeler\thanks{Partially supported by the European Commission
under grant HPRN-CT-1999-00118}, \, J.~Latschev$^\dag$ and
E.~Zehnder}
\date{February 10, 2003}
\begin{document}
\maketitle
\sloppy
\begin{abstract} In this paper we find conditions which guarantee that a given flow $\Phi$ on a compact metric space $X$ 
admits a Lyapunov one-form $\omega$ lying in a prescribed \v Cech cohomology class $\xi\in \check H^1(X;\R)$. These conditions are formulated in terms
of the restriction of $\xi$ to the chain recurrent set of $\Phi$. The result of the paper may be viewed as a generalization of a well-known theorem
of C. Conley about the existence of Lyapunov functions. 
\end{abstract}

\section{Introduction}
\label{intro}

C. Conley proved in \cite{C, C1} that any flow $\Phi: X\times \R
\to X$ on a compact metric space $X$ decomposes into a chain
recurrent flow and a gradient-like flow. More precisely, he
proved the existence of a {\it Lyapunov function} for the flow, i.e. a
continuous function $L: X\to \R$, which decreases along any orbit
of the flow lying in the complement $X-R$ of the chain recurrent set $R\subset X$ of $\Phi$ and is constant on the connected 
components\footnote{C. Conley (see \cite{C1}, Theorem 3.6D) proved that the chain transitive components of $R$ 
coincide with the connected
components of $R$.}
of $R$.

\begin{theorem} {\bf (C. Conley, {\rm \cite{C,C1}})}
\label{intro1} Let $\Phi:X \x \bbr \to X$, $\Phi(x,t)=x\cdot t$, be a continuous flow on
a compact metric space $X$. Then there exists a continuous
function $L:X\to \bbr$, which is constant on the connected
components of the chain recurrent set $R=R(\Phi)$ of the flow
$\Phi$, and satisfies $ \label{introeq.1} L(x\cdot t) <L(x)$ for
any $x\in X - R$ and $t>0$. 
\end{theorem}
This important result led Conley to his program of understanding
very general flows as collections of isolated invariant sets
linked by heteroclinic orbits.

Our aim in this paper is to go one step further and to analyze the
flow within the chain-recurrent set $R$, where it is typically
complicated. As a new tool, we study the notion of a
{\it Lyapunov one-form} for $\Phi$, which is a natural
generalization of the notion of a Lyapunov function and has been
introduced in a different context in Farber's papers \cite{Fa:1, Fa:2}. We prove
that under some natural assumptions, a given \v Cech cohomology
class $\xi\in \check H^1(X;\R)$ can be represented by a continuous
closed Lyapunov one-form for the flow $\Phi$.

The notion of a {\it continuous closed one-form on a topological
space} generalizes the notion of a continuous function. For
convenience of the reader we recall the relevant
definitions and the main properties of continuous closed one-forms
in \S \ref{prelim}, referring for more details to the papers \cite{Fa:1, Fa:2}, 
where they were originally introduced.
For the purpose of this
introduction, let us say that continuous closed one-forms
are analogues of the familiar smooth closed one-forms on
differentiable manifolds. Any continuous closed one-form $\omega$
on a topological space $X$ canonically determines a \v Cech
cohomology class $[\omega]\in \check H^1(X;\R)$, which plays a
role analogous to the de Rham cohomology class of a smooth closed
one-form. For any continuous curve $\gamma:[0,1]\to X$, the line
integral $\int_\gamma \omega \in \R$ is defined and has the usual
properties; in particular, it depends only on the homotopy class
of the curve relative to the end-points.

\begin{definition}
\label{introlyap}
Consider a continuous flow $\Phi: X\times \R\to X$ on a topological space $X$.
Let $Y\subset X$ be a closed subset invariant under $\Phi$. A
continuous closed one--form $\omega$ on $X$ is called a Lyapunov one--form for the pair
$(\Phi, Y)$ if it has the following two properties:
\begin{enumerate}
\item[{\rm (L1)}] For every $x\in X- Y$ and every $t>0$,
$$
\label{introeq.2} \int_x^{x\cdot t} \omega \, <\, 0,
$$
where the integral is calculated along the trajectory of the flow.
\item[{\rm (L2)}] There exists a continuous function 
$f:U\to \R$ defined on an open neighborhood $U$ of $Y$ such that 
$\omega|_U =df$ 
and $f$ is constant on any connected component of $Y$.
\end{enumerate}
\end{definition}

Any continuous function $L:X\to \R$ determines the closed one-form
$\omega=dL$ (see \S \ref{prelim}) and in this special case, condition (L1)
reduces to the requirement $L(x\cdot t)<
L(x)$ for any $t>0$ and $x\in X- Y$, while condition (L2) means 
that $L$ is constant on any connected component of $Y$.
Hence for $\omega=dL$,
the above definition reduces to
the classical notion of a Lyapunov function, see \cite{Sh}.

The following remark illustrates Definition \ref{introlyap}. Given
a flow $\Phi$ on $X$ and a Lyapunov one--form $\omega$ for
$(\Phi,Y)$, representing a nonzero \v Cech cohomology class
$[\omega] =\xi \in \check H^1(X;\bbr)$, the homology class $z\in
H_1(X; \bbz)$ of any periodic orbit of $\Phi$ satisfies
$$ \langle \xi, z \rangle \leq 0,$$
with equality if
and only if the periodic orbit is contained in $Y$. Using this fact one constructs flows 
such that no nonzero cohomology class $\xi\in \check H^1(X;\bbr)$ contains a Lyapunov 1-form.

In this paper we will associate with any cohomology class $\xi\in
\check H^1(X;\R)$ a subset $R_\xi\subset R$ of the chain recurrent
set $R=R(\Phi)$ of the flow $\Phi$, see section \S \ref{unfold}
for details. The set $R_\xi$ is closed and invariant under the flow and can be characterized as the projection of
the chain recurrent set of the natural lift of the flow to the
Abelian cover of $X$ associated with the class $\xi$. 

\begin{figure}[h]
\begin{center}
\includegraphics[0,0][215,130]{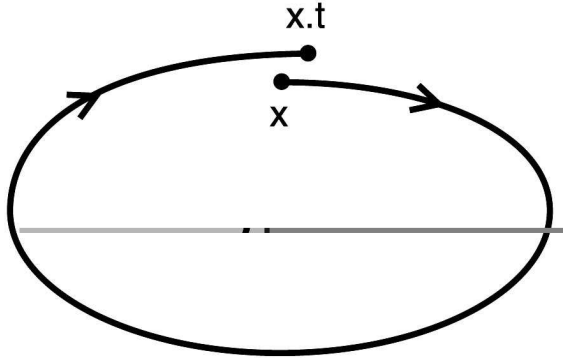}
\end{center}
\caption{$(\delta,T)$-cycle}\label{cycle}
\end{figure}

A {\it $(\delta,T)$-cycle} is a pair $(x,t)\in X \x \bbr$ satisfying
$t\geq T$ and $d(x,x\cdot t) <\delta$.  Here $d$ denotes the
distance function on $X$. See Figure \ref{cycle}. If $X$ is
locally path-connected, any $(\delta, T)$-cycle with small enough
$\delta$ determines a closed loop, which first follows the flow
line from $x$ to $x\cdot t$ and then returns from $x \cdot t$ to
$x$ by a path contained in a suitably small ball. This leads to
the notion of {\it a homology class $z\in H_1(X;\Z)$ associated to
a $(\delta,T)$-cycle}, see Definition \ref{associated}. The class
$z$ is uniquely defined if $X$ is homologically locally
1-connected; without this assumption the homology class $z$
associated with a $(\delta, T)$-cycle might not be unique.
The natural bilinear pairing
$\langle \,  , \,  \rangle : \check H^1(X; \bbr) \x H_1(X; \bbz) \to\bbr$
can be understood as
$
\langle \xi, z \rangle = \int\limits_\gamma \omega,
$
where $\omega$ is a representative closed one--form for $\xi \in \check
H^1(X; \bbr)$ and $\gamma:[0,1] \to X$ is a loop representing the class $z\in
H_1(X; \bbz)$.
Despite the fact that the homology class $z$ associated to a $(\delta,T)$-cycle might depend
(for wild $X$) on the choice of the connecting path between $x \cdot
t$ and $x$, the construction is such that the value $\langle \xi, z \rangle
\in \bbr$ only depends (for small enough $\delta>0$) on the $(\delta,T)$-cycle itself; see \S
\ref{unfold} for details.

The following Theorem is our main result:

\begin{theorem}
\label{intromain}
Let $\Phi$ be a continuous flow on a compact, locally path connected, metric space
$X$ and $\xi$ a cohomology class in $\check H^1(X; \bbr)$. 
Denote by $C_\xi$ the subset 
\begin{eqnarray}
C_\xi = R - R_\xi\label{circulation}
\end{eqnarray}
of the chain recurrent set $R$ of the flow $\Phi$.
Assume that the following two conditions are satisfied
\begin{enumerate}[{\rm (A)}]
\item $\xi|_{R_\xi} = 0$;
\item there exist constants $\delta>0$, $T>1$, such that every
homology class $z\in H_1(X; \bbz)$ associated to an arbitrary
$(\delta,T)$-cycle $(x,t)$ with $x\in C_\xi$ satisfies
$$ \langle \xi , z
\rangle \leq -1.
$$
\end{enumerate}
Then there exists a Lyapunov one-form $\omega$ for $(\Phi,R_\xi)$ representing
the cohomology class  $\xi$. 
Moreover, the subset $C_\xi$ is closed.

Conversely, if for the given cohomology class $\xi$ 
there exists a Lyapunov one-form for the pair $(\Phi,R_\xi)$ in the class $\xi$ 
and if the set $C_\xi$ is closed 
then {\rm (A)} and {\rm (B)} hold true.
\end{theorem}

\begin{corollary}\label{removed}
Let $\Phi: X\times \R\to X$ be a continuous flow on a compact,
locally path connected metric space. Any \v Cech cohomology
class $\xi\in \check H^1(X;\R)$ satisfying $\xi|_R=0$ (where $R
=R(\Phi)$ denotes the chain recurrent set of the flow), contains
a Lyapunov one-form $\omega$ for $(\Phi, R)$.
\end{corollary}

Corollary \ref{removed} follows directly from Theorem
\ref{intromain} since under the assumption $\xi|_R=0$ the set
$R_\xi$ coincides with $R$ (compare Definition \ref{def4.2}) and
so the set $C_\xi$ is empty. Corollary \ref{removed} also admits a
simple, independent proof based on Conley's Theorem \ref{intro1}.

\begin{corollary}
\label{introcor1} Suppose $\Phi:X \x \bbr \to X$ is a flow on a
compact locally path connected metric space, whose chain recurrent
set consists of finitely many rest points and periodic orbits.
Then a Lyapunov one--form  for $(\Phi, Y)$, with suitable
$Y\subset X$, exists in a nontrivial cohomology class $\xi\in
\check H^1(X; \bbr)$ if and only if the homology classes of the
periodic orbits are contained in the half space
$$
H_\xi:= \{ z \in H_1(X; \bbz) \; | \; \langle \xi, z \rangle \leq 0 \}.
$$
If the above condition holds, then the set $Y$ coincides with the
union of the rest points and of those periodic orbits, for which
the corresponding homology classes $z\in H_1(X;\Z)$ satisfy
$\langle \xi, z\rangle =0$.
\end{corollary}

This Corollary is a direct consequence of our main Theorem
\ref{intromain}.
The class of flows meeting its assumptions 
includes the Morse-Smale flows on closed manifolds.

Another interesting special case arises when $R_\xi=\emptyset$.
From Theorem \ref{intromain} we also deduce the following result: 

\begin{corollary}\label{introcor.2}
Let $\Phi: X\times \R \to X$ be a continuous flow on a compact,
locally path connected, metric space $X$ and let $\xi\in \check
H^1(X; \R)$ be a nonzero \v Cech cohomology class. The following
two conditions are equivalent:
\begin{enumerate}
\item[\rm {(i)}] $R_\xi=\emptyset$ and the flow satisfies condition $(B)$ of Theorem
\ref{intromain};
\item[{\rm (ii)}] there exists a Lyapunov one--form for $(\Phi,\emptyset)$
representing the class $\xi$.
\end{enumerate}
If the class $\xi$ is integral, i.e. $\xi\in \check
H^1(X;\Z)$, then either of the above conditions is equivalent to
the existence of a continuous locally trivial fibration $p:X \to
S^1\subset \C$ with the following properties. The function $t\mapsto \arg (p(x\cdot
t))$ is differentiable, the derivative
$$\displaystyle{\frac{d}{dt}} \arg (p(x \cdot t)) <0$$
is negative for all $x\in X$, $t\in \R$, and the cohomology class
 $p^*(\mu)\in \check
H^1(X;\Z)$ coincides with $\xi$, where $\mu\in \check H^1(S^1;
\Z)$ is the fundamental class of the circle $S^1$ oriented
counterclockwise. In particular, for any angle $\theta\in S^1$ the
set $K=p^{-1}(\theta)$ is a cross section of the flow $\Phi$.
\end{corollary}

Recall that a closed subset $K \subset X$ is a {\it cross section} of
the flow $\Phi$ if the flow map $K \x \bbr \to X$ is a surjective
local homeomorphism (see S. Schwartzman \cite{Sc}). A cross section $K$ is transversal to the flow and the orbit of 
every
point in $X$ intersects $K$ in forward and backward time.

We give a proof of Corollary
\ref{introcor.2} in section \ref{applic}.

Under slightly stronger assumptions, S. Schwartzman \cite{Sc1} proved, among other results, the equivalence
of the property of having a cross section and the existence of a fibration $p:X\to S^1$
as stated in the second part of Corollary \ref{introcor.2}.
In subsequent work we will compare in detail our results with Schwartzman's beautiful paper \cite{Sc}.

The second part of Corollary \ref{introcor.2} may be viewed as a generalization of
Fried's results on the existence of cross sections to flows on
manifolds \cite[Theorem D]{Fr:1}. 
The assumption of D. Fried to ensure the existence of cross sections
is formulated in terms of the notion of {\it homological directions of a flow}, which we now recall. 
A sequence $(x_n, t_n)\in X\times (1, \infty)$ is a {\it closing sequence based at $x\in X$} if the sequences 
$x_n$ and $x_n\cdot t_n$ tend to $x$. We will assume that $X$ is a compact polyhedron.
Then any closing sequence determines  uniquely a sequence
of homology classes $z_n\in H_1(X;\Z)$. Here $z_n$ denotes the homology class of a loop,
which starts at $x_n$, follows the flow until 
$x_n\cdot t_n$, and then returns to $x_n$ along a "short" path. 
Let $D_X$ be the factorspace $D_X= H_1(X;\R)/\R_+$, where this space is topologized as the disjoint union 
of the unit sphere with the origin. Any closing sequence as above determines a sequence of "homology directions" 
$\tilde z_n\in D_X$, the equivalence classes of $z_n$ in $D_X$.
The set of homology directions $D_\Phi\subset D_X$ of the flow $\Phi$ is 
defined as the set of all accumulation points of all sequences $\tilde z_n$ corresponding to closing sequences in $X$. 
As noted by D. Fried \cite{Fr:1} it is enough to consider closing sequences $(x_n,t_n)$ with $t_n\to \infty$. 

\begin{prop}\label{closing}
Let $X$ be a finite polyhedron and let $\xi\in H^1(X;\Z)$ be an integral cohomology class.
 Let $\Phi: X\times \R \to X$ be a continuous flow such that the chain recurrent set $R_\xi$ is isolated in $R$. 
Then condition (B) of Theorem 
\ref{intromain} is equivalent to the Fried's condition that any homology direction $\tilde z= \lim \tilde z_n\in D_X$
of any closing sequence 
$(x_n, t_n)\in X\times (1, \infty)$ with $x_n \in C_\xi$ satisfies 
$\langle \xi, \tilde z\rangle <0.$
\end{prop}

See \S \ref{proofclosing} for a proof.

A comparison of the results of D. Fried \cite{Fr:1} with the results of this paper shows that
our setting is more general in two respects: we allow spaces $X$ of a more general nature and arbitrary real \v Cech cohomology classes
$\xi$. In \cite{Fr:1} $X$ is required to be a compact manifold, possibly with boundary, and the class $\xi$ has to be integral.
The equivalence of our condition (B) with the condition of D. Fried (Proposition \ref{closing})
holds only under these additional assumptions.

In the papers \cite{Fa:1}, \cite{Fa:2} two different generalizations ${\rm cat}(X,\xi)$ and ${\rm Cat}(X,\xi)$ of the classical notion of Lusternik - Schnirelman 
category ${\rm cat}(X)$ were introduced; here $X$ is a finite polyhedron and $\xi\in H^1(X;\R)$ a cohomology class. 
Using these new concepts, an extension of the Lusternik - Schnirelman theory for flows was constructed, see \cite{Fa:1}, \cite{Fa:2}. 
The main results of \cite{Fa:1}, \cite{Fa:2} allow to estimate the number of fixed points of a flow under the assumption that
(1) the fixed points are isolated in the chain recurrent set, and (2) the flow admits a Lyapunov closed 1-form 
lying in the class $\xi$. The results of the present paper fit nicely in this program and explain the nature of the assumption (2).
Note that (1) is similar in spirit (although formally not equivalent) to the property that the set $R_\xi$
is isolated in the chain recurrent set $R$ of the flow.

\section{Closed one-forms on topological spaces}
\label{prelim}

In this section we recall the notion of a continuous closed one-form on a topological space,
which has been introduced in \cite{Fa:1, Fa:2}.

\begin{definition} Let $X$ be a topological space.
A {\it continuous closed one-form} on $X$ is defined by a an open cover
$\cu = \{ U \}$ of $X$ and by a collection
$\{\pHi_U\}_{U\in \cu}$ of continuous functions $\pHi_U:U \to \bbr$ with the following
property: for any two subsets $U,V \in \cu$ the difference
\begin{eqnarray}
\pHi_{U}|_{U \cap V} - \pHi_{V}|_{U\cap V} : U \cap V \to \bbr\label{locallyconst}
\end{eqnarray}
is a locally constant function (i.e. constant on each connected component of $U\cap V$).
Two such collections $\{\pHi_U\}_{U\in \cu}$ and $\{\psi_V\}_{V\in \mathcal V}$
are called equivalent if their union $\{\pHi_U, \psi_V\}_{U\in \cu, V\in \mathcal V}$
satisfies condition (\ref{locallyconst}). The equivalence classes are
called continuous closed one--forms on $X$.
\end{definition}

Any continuous function $f:X \to \bbr$ (viewed as a family consisting of a single element)
determines a continuous
closed one--form, which we denote by $df$ (the differential of $f$).

A continuous closed
one--form $\omega$ vanishes, $\omega=0$, if it is represented by a collection of locally
constant functions.

The sum of two continuous closed one--forms determined by collections
$\{\pHi_U\}_{U\in\cu}$ and $\{\psi_V\}_{V\in \cv}$ is the continuous closed
one--form corresponding to the collection $\{\pHi_{U|U \cap V} + \psi_{V|U
\cap V} \}_{U\in \cu, V\in \cv}$. Similarly, one can multiply continuous
closed one--forms by real numbers $\lambda \in \bbr$ by multiplying the
corresponding representatives with $\lambda$. With these operations, the set
of continuous closed one--forms on $X$ is a real vector space.

Continuous closed one--forms behave naturally with respect to continuous maps:
If $h:Y \to X$ is continuous and $\{\pHi_U\}_{U\in \cu}$ determines a
continuous closed one--form $\omega$ on $X$, then the collection
of continuous functions $\phi_U\circ h: h^{-1}(U)\to \R$ determines a
continuous closed one--form on $Y$ which will be denoted
by $h^*\omega$. As a special case of this construction, we will often use the operation of
restriction of a closed one-form $\omega$ to a given subset
$A\subset X$; in this case $h$ is the inclusion map $A\to X$ and the form $h^\ast\omega$ is simply denoted as
$\omega|_A$.

A continuous closed one--form $\omega$ on $X$ can be integrated
along continuous paths in $X$. Namely, let $\omega$ be given by a
collection $\{\pHi_U\}_{U\in \cu}$, and 
$\gamma:[0,1] \to X$ be a continuous path. We may find a finite subdivision
$0=t_0<t_1 < \dots < t_N = 1$ of the interval $[0,1]$ such that
for each $1\leq i \leq N$ the image $\gamma([t_{i-1},t_{i}])$ is
contained in a single open set $U_i\in \cu$. Then we define the
line integral
\begin{equation}
\label{eq1.2} \int_\gamma \omega := \sum_{i=1}^{N}
\pHi_{U_i}(\gamma(t_{i})) - \pHi_{U_i}(\gamma(t_{i-1})).
\end{equation}
The standard arguments show that the integral (\ref{eq1.2}) is independent
of all choices and in fact depends only on the homotopy class of $\gamma$
relative to its endpoints.

Consider the following exact sequence of sheaves over $X$
\begin{eqnarray}
0\to \R_X \to C_X \to B_X\to 0,\label{sheaves}
\end{eqnarray}
where $\R_X$ is the sheaf of locally constant functions, $C_X$ is the sheaf of real-valued continuous functions,
and $B_X$ is the sheaf of germs of continuous functions modulo locally constant ones. More precisely,
$B_X$ is the sheaf corresponding to the presheaf $U\mapsto C_X(U)/\R_X(U)$. By the definitions above, 
the global sections of the sheaf $B_X$
are in one-to-one correspondence with continuous closed 1-forms on $X$. Hence the space of all closed 1-forms on $X$ is $H^0(X;B_X)$.

The exact sequence of sheaves (\ref{sheaves}) generates the cohomological exact sequence
\begin{eqnarray}
0\to H^0(X;\R_X) \to H^0(X;C_X)\stackrel d\to H^0(X;B_X)\stackrel {[\, \, ]}\to H^1(X;\R_X)\to 0\label{sheaf1}
\end{eqnarray}
In this exact sequence, $H^0(X;C_X)=C(X)$ is the space of all continuous functions $f: X\to \R$ and
the map $d$ assigns to any continuous function $f$ its differential $df\in H^0(X;B_X)$. The group
$H^1(X;\R_X)$ is the \v Cech cohomology $\check H^1(X;\R)$ (see \cite{Sp}, chapter 6);
the map $[\, \, ]$ assigns to any closed
one-form $\omega$ its \v Cech cohomology class $[\omega]\in \check H^1(X;\R)$.
This proves:

\begin{lemma}\label{check}
A continuous closed 1-form $\omega\in H^0(X;B_X)$ equals $df$ for some continuous function $f:X\to \R$
if and only if its \v Cech cohomology class $[\omega]\in \check H^1(X;\R)$ vanishes, $[\omega]=0$.
Any \v Cech cohomology class $\xi\in \check H^1(X;\R)$ can be realized by a continuous closed 1-form on $X$.
\end{lemma}

Note also that there is a natural homomorphism $\check H^1(X;\R)\to H^1(X;\R)$ from \v Cech cohomology to
singular cohomology. Using Lemma \ref{check} and the well-known identification $H^1(X;\R)\simeq \Hom(H_1(X);\R)$, it can be described as a pairing
\begin{eqnarray}
\langle \, ,\, \rangle : \check H^1(X;\R)\times H_1(X;\Z)\to \R, \quad \mbox{where}\quad
\langle [\omega], [\gamma]\rangle = \int_\gamma \omega.
\end{eqnarray}
In other words, choosing a representative closed 1-form $\omega$ for a cohomology class
$\xi\in \check H^1(X;\R)$ and a closed loop $\gamma$ in $X$
representing a homology class $z\in H_1(X;\Z)$, the number $\langle \xi, z\rangle \in \R$ equals the line integral
$\int_\gamma \omega$, which is independent of the choices.

Here is a generalization of the well-known {\it Tietze
extension theorem}:

\begin{prop}
\label{tietze}
Let $X$ be a metric space and $A\subset X$ a closed subset. Let
$\omega$ be a continuous closed one--form on $A$, and let $\xi\in \check
H^1(A;\R)$ denote the \v Cech cohomology class of $\omega$. Then for any
cohomology class $\xi'\in \check H^1(X;\R)$, satisfying $\xi'|_A =\xi$, there
exists a continuous closed one--form $\omega'$ on $X$ representing the
cohomology class $\xi'$, such that $\omega'|_A=\omega$.
\end{prop}

\begin{proof} Choose an arbitrary continuous closed one--form $\Omega'$
representing the class $\xi'$. Then $\Omega'|_A$ is cohomologous to $\omega$,
i.e. $\Omega'|_A - \omega =df$, where $f: A\to \R$ is a continuous
function. By the Tietze Extension Theorem for functions, we find a
continuous function $f': X\to \R$ extending $f$. Then $\omega' = \Omega'
-df'$ is a closed one--form in the class $\xi'$ satisfying $\omega'|_A=\omega$.
$\square$
\end{proof}

The statement of Proposition \ref{tietze} can be expressed as follows: a continuous closed 1-form
$\omega$ on a closed subset $A\subset X$ can be extended to a continuous closed 1-form on $X$ if and only if
the cohomology class $[\omega]\in \check H^1(A;\R)$ can be extended to a cohomology class lying in
$\check H^1(X;\R)$.

\section{The chain recurrent set $R_\xi$}
\label{unfold}

The goal of this section is to introduce a new chain recurrent set
$R_\xi=R_\xi(\Phi)\subset X$, which is associated with a flow
$\Phi$ on $X$ together with a \v Cech cohomology class $\xi\in \check
H^1(X;\R)$. The set $R_\xi$ appears in the statement of our main
Theorem \ref{intromain}.

Throughout this section we will assume that $X$ is a locally path connected
compact metric space.

Recall that a space $X$ is {\em locally path connected} if
for every open set $U\subset X$ and for every point $x\in U$ there exists an open set
$V\subset U$ with $x\in V$ such that any two points in $V$ can be connected
by a path in $U$. Equivalently, $X$ is locally path connected iff the connected components of open subsets
are open (see \cite{Sp}, page 65).

\subsection{Definition of $R_\xi$}\label{defxi}

Recall the definition of the chain-recurrent set $R=R(\Phi)$ of
the flow $\Phi$. Given any $\delta>0$, $T>0$, a $(\delta,T)$-{\it chain from $x\in X$ to $y\in
X$}  is a finite sequence $x_0=x, x_1, \dots, x_{N}=y$ of points in $X$
and numbers $t_1, \dots, t_N\in \R$ satisfying $t_i\geq T$ and
$d(x_{i-1}\cdot t_i, x_i) < \delta$ for all $1 \leq i \leq N$.
\begin{figure}
\begin{center}
\includegraphics[0,0][344,64]{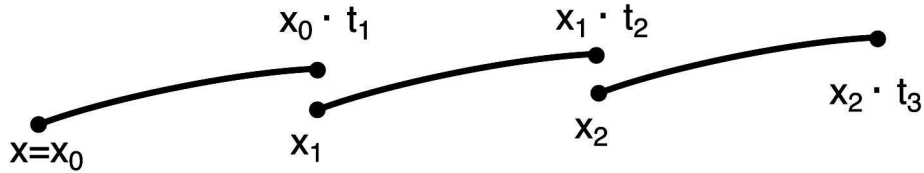}
\end{center}
\caption{$(\delta,T)$-chain}\label{chain}
\end{figure}
Note that a $(\delta,T)$-cycle (see \S \ref{intro}) is a $(\delta,
T)$-chain of a special kind (with $n=1$).

The {\it chain recurrent set} $R=R(\Phi)$ of the flow
$\Phi$ is defined as the set of all points $x\in X$ such that for any
$\delta>0$ and $T \geq 1$ there exists a $(\delta,T)$-chain starting and ending
at $x$. It is immediate from this definition that the chain recurrent set is
closed and invariant under the flow; $R$ contains all
fixed points and periodic orbits.
The chain recurrent set $R$
contains the set of all nonwandering points and, in particular, the positive
and negative limit sets of any orbit \cite[\S II.6]{C}.
The set $R=R(\Phi)$ is a disjoint union of its {\it chain transitive components}\footnote
{Recall that $x, y\in R$ belong to the same chain transitive component if for any $\delta>0$ and $T>1$ there exist
$(\delta, T)$-chains from $x$ to $y$ and from $y$ to $x$.}.

\begin{lemma}
\label{lem4.1}
Given a locally path connected compact metric space $X$ and
a number $\eps>0$, there exists $\delta = \delta(\eps)>0$ such that
any two points $x, y\in X$ with $d(x,y)<\delta$ can be connected by a continuous path
$\gamma: [0,1]\to X$ contained in some open $\eps$-ball.
\end{lemma}

\begin{proof}
By the definition of local path-connectedness, each point $x\in X$ has a
neighborhood $V_x$ contained in the $\eps$-ball $B_x$ around $x$ such that
any two points in $V_x$ can be connected by a path in $B_x$. Choose a finite
subcover of the covering $\{V_x\}_{x\in X}$ of $X$ and choose 
$\delta(\eps)$ as the Lebesgue number of this finite cover. $\square$
\end{proof}

\begin{definition}\label{scale}
A pair $(\varepsilon, \delta)$ 
of real numbers $\varepsilon =\varepsilon(\xi)>0$ and $\delta=\delta(\xi)>0$ is called a scale of a nonzero cohomology class 
$\xi\in \check H^1(X;\R)$ if {\rm (1)} $\xi|_{B}=0$ for any ball $B\subset X$ of radius $2\varepsilon$ and {\rm (2)}
any
two points $x, y\in X$ with $d(x,y)<\delta$ can be connected by a path in $X$
contained in a ball of radius $\eps$.
\end{definition}

Such a scale always exists. In fact, since we may realize the class $\xi$
by a continuous closed 1-form
$\omega=\{\phi_U\}_{U\in \U}$ with a finite open cover $\U$ and then take for $\varepsilon$ half of the Lesbegue
number of $\U$. 
Using Lemma \ref{lem4.1} we may then find $\delta =\delta(\xi)>0$ satisfying condition (2) in Definition \ref{scale}.

We want to evaluate \v Cech cohomology classes on
broken chains of trajectories of the flow which start and end at the same point. This can be done as follows.
Let $\varepsilon=\varepsilon(\xi)$, $\delta=\delta(\xi)$ be a scale for $\xi\in \check H^1(X;\R)$ (see Definition \ref{scale}).
Suppose we are given a {\it closed $(\delta,T)$-chain},
i.e. a $(\delta,T)$-chain from a point $x$ to itself.
We have a sequence of points $x_0=x, x_1, \dots, x_{N-1}, x_{N}=x$ of $X$
and a sequence of numbers $t_1, \dots, t_N\in \R$ with $t_i\geq T$, such that $d(x_{i-1}\cdot t_i, x_i) < \delta$ for
any $1 \leq i \leq N$. We want to associate with such a chain a homology class
$z \in H_1(X; \bbz)$.
Choose continuous paths $\sigma_i:[0,1]\to X$, where $1 \leq i \leq N$, connecting
$x_{i-1}\cdot t_i$ with $x_i$ and lying in a ball $B_i$ of radius $\eps$. We obtain a singular cycle
which is a combination of the parts of the trajectories from $x_{i-1}$ to $x_{i-1}\cdot t_i$ and the paths $\sigma_i$.
\begin{definition}\label{associated}
The homology class $z \in H_1(X; \bbz)$
of this singular cycle is said to be associated with the given closed $(\delta, T)$-chain.
\end{definition}

Note that the obtained class $z$ may depend on the choice of paths $\sigma_i$ (if the space $X$ is wild, i.e. not locally contractible). However the value
\begin{equation}
\label{eq4.1}
\langle \xi, z\rangle \in \R, \quad\mbox{where}\quad  \langle \xi, z\rangle =  \sum_{i=1}^{N}
\int_{x_{i-1}}^{x_{i-1} \cdot t_i} \omega \quad + \quad
\sum_{i=1}^N \int_{\sigma_i} \omega,
\end{equation}
is independent of the paths $\sigma_i$. Indeed, if we use two different sets of curves $\sigma_i$ and  $\sigma_i'$
then the difference of the corresponding expressions in (\ref{eq4.1}) will be the integral over the sum of singular cycles
$\sum_{i=1}^N (\sigma_i -\sigma_i')$, each cycle $\sigma_i-\sigma_i'$ being contained in a ball of radius
$2\varepsilon$. Since we know that the restriction of the cohomology class $\xi$ on any such ball vanishes, we see that the right side of
(\ref{eq4.1})
is independent of the choice of the curves $\sigma_1, \dots, \sigma_N$.

The homology class $z\in H_1(X;\Z)$ associated with a closed  $(\delta,T)$-chain
is uniquely defined if $X$ is homologically locally connected in dimension 1 (see \cite{Sp}, chapter 6, page 340 for the definition)
and $\delta>0$ is sufficiently small.

Now we are ready to define the subset $R_\xi$ of the chain reccurent set $R$:
\begin{definition}
\label{def4.2}
Let $\varepsilon=\varepsilon(\xi)$ and $\delta=\delta(\xi)$ be a scale of a cohomology class $\xi\in \check H^1(X;\R)$.
Then $R_\xi = R_\xi(\Phi)$ denotes the set of all points $x\in X$ with
the following property: for any $0<\delta'<\delta$ and $T>1$ there exists a
$(\delta',T)$-chain from $x$ to $x$ such that $\langle \xi, z\rangle = 0$
for any  homology
class $z\in H_1(X; \bbz)$ associated with this chain.
\end{definition}

Roughly, the set $R_\xi$ can be characterized as the part of the chain recurrent set of the flow
in which the cohomology class $\xi$ does not detect the motion.

It is easy to see that $R_\xi$ is closed and invariant with respect to the flow.

Note also that $R_\xi =R_{\xi'}$ whenever $\xi'
=\lambda\xi$, with $\lambda\in \R$, $\lambda\neq 0$. Thus, $R_\xi$ depends
only on the line through $\xi$ in the real vector space $\check H^1(X;\R)$.

Any fixed point of the flow belongs to $R_\xi$. The points of a
periodic orbit belong to $R_\xi$ if the homology class
$z\in H_1(X;\Z)$ of this orbit satisfies $\langle \xi, z\rangle
=0$.

{\bf Example.}
\begin{figure}[h]
\begin{center}
\includegraphics[121,-313][251,-183]{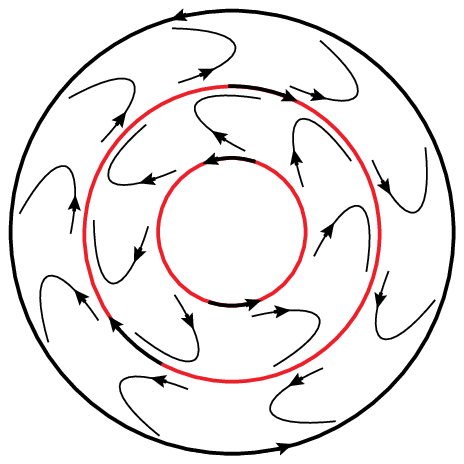}
\end{center}
\caption{The flow on the planar ring $Y$.}\label{example2}
\end{figure}
It may happen that the points of a periodic orbit belong to $R_\xi$ although 
$\langle \xi, z\rangle\not=0$ for the homology class $z$
of the orbit. This possibility is illustrated by the following example.  
 
Consider the flow on the planar ring $Y\subset \C$ shown on Figure \ref{example2}.
In the polar coordinates $(r, \phi)$ the ring $Y$ is given by the inequalities $1\leq r\leq 3$ and the flow is given by the differential equations
$$\dot r = (r-1)^2(r-3)^2(r-5)^2, \quad \dot\phi =\sin\left(r\cdot \frac{\pi}{2}\right).$$
Let $C_k$, where $k=1, 2, 3$, denote the circle $r=2k-1$.
The circles $C_1, C_2, C_3$ are invariant under the flow. 
The motion along the circles $C_1$ and $C_3$ has constant angular velocity 1.   
Identifying any point $(r, \phi)\in C_1$ with $(5r, \phi)\in C_3$
we obtain a torus $X=Y/\simeq$ and a flow $\Phi: X\times \R\to X$. 
The images of the circles $C_1, C_2, C_3\subset Y$ represent two circles $C'_1=C'_3$ and $C'_2$ on the torus $X$. 

Let $\xi\in H^1(X;\R)$ be a nonzero cohomology class which is the pullback of a cohomology class of $Y$.
One verifies that in this example the set $R_\xi(\Phi)$ coincides with the whole torus $X$. 
In particular, $R_\xi(\Phi)$ contains the periodic orbits $C_1'=C'_3$ and $C_2'$ although clearly
$\langle \xi, [C'_k]\rangle \not=0$.

\subsection{$R_\xi$ and dynamics in the free Abelian cover}\label{dynamics}

Now we will give a different characterization of $R_\xi$ using the
dynamics in the covering space associated with the class $\xi$. We
will use chapter 2 in \cite{Sp} as a reference for notions
related to the theory of covering spaces.

Recall our standing assumption that $X$ is a locally path
connected compact metric space. For simplicity of exposition we will additionally assume that $X$ is connected.

Any \v Cech cohomology class $\xi\in \check H^1(X;\R)$ determines a homomorphism
\begin{eqnarray}
h_\xi : \pi_1(X,x_0) \to \R, \quad h_\xi([\alpha]) = \langle \xi, [\alpha]\rangle =\int_\alpha \omega \, \in \R,\label{periods}
\end{eqnarray}
where $\alpha :[0,1]\to X$ is a continuous loop $\alpha(0)=\alpha(1)=x_0$,
$[\alpha]\in \pi_1(X, x_0)$ denotes its homotopy class, and $\omega$ is a continuous closed one-form in the class $\xi$.
The map $h_\xi$ is called {\it homomorphism of periods}.

The kernel of $h_\xi$ is a normal subgroup $H=\ker (h_\xi) \subset
\pi_1(X,x_0)$. We want to construct a covering projection map $p_\xi:
\tilde X_\xi \to X$, corresponding to $H$, i.e. $(p_\xi)_{\#}\pi_1(\tilde X_\xi, \tilde
x_0) =H$. The uniqueness of such a covering projection map follows from
\cite{Sp}, chapter 2, Corollary 3. To show the existence we may use
\cite{Sp}, chapter 2, Theorem 13; according to this Theorem (see
also Lemma 11 in \cite{Sp}, chapter 2) we have to show that for
some open cover $\U$ of $X$ the subgroup $\pi_1(\U, x_0)\subset
\pi_1(X,x_0)$ is contained in $H$. Here $\pi_1(\U, x_0)\subset
\pi_1(X,x_0)$ denotes the subgroup generated by homotopy classes
of the loops of the form $\alpha = (\gamma\ast\gamma')\ast
\gamma^{-1}$ where $\gamma'$ is a closed loop lying in some
element of $\U$ and $\gamma$ is a path from $x_0$ to $\gamma'(0)$.
To show that this condition really holds in our situation, let us
realize $\xi$ by a closed one-form $\omega=\{f_U\}_{U\in \U}$, $\U$ being an open cover of $X$. 
We claim that $\pi_1(\U, x_0)\subset H$ for this cover $\U$. Indeed, for any loop of the form $\alpha
= (\gamma\ast\gamma')\ast \gamma^{-1}$ where $\gamma'$ lies in
some $U\in \U$,
$$\langle \xi, [\alpha]\rangle =\int_\alpha \omega = \int_{\gamma'}\omega =0,$$
since $\gamma'$ lies in $U$ and $\xi|_U=0$.
Thus:
{\it Any \v Cech cohomology class $\xi\in \check H^1(X;\R)$ uniquely determines a covering projection
map $p: \tilde X_\xi \to X$ with connected total space $\tilde X_\xi$, such that 
$(p_\xi)_{\#}\pi_1(\tilde X_\xi, \tilde x_0) =\ker (h_\xi)$.}

\begin{lemma} Let $X$ be a connected and locally path connected compact metric space
and let $\xi\in \check H^1(X;\R)$ be a \v Cech cohomology class. The group of covering transformations of the covering
map $\tilde X_\xi \to X$ is a finitely generated free Abelian group.
\end{lemma}
\begin{proof}
The group of covering transformations of $\tilde X_\xi \to X$ can be identified with
$\pi_1(X,x_0)/H$, which is isomorphic to the image
of the homomorphism of periods $h_\xi(\pi_1(X,x_0))\subset \R.$ It
is a subgroup of $\R$ and hence it is Abelian and has no torsion.
Therefore it is enough to show that the image of the homomorphism
of periods $h_\xi(\pi_1(X,x_0))\subset \R$ is finitely generated.

Let $\omega=\{f_U\}_{U\in \U}$ be a continuous closed 1-form with respect to an open cover $\U$ representing $\xi$.
Find an open cover $\V$ of $X$ and a function $\kappa : \V
\to \U$ such that for any $V\in \V$, the set $U=\kappa(V) \in \U$
satisfies $\bar V\subset U$. We may realize $\omega$ with respect
to the open cover $\V$ as $\omega=\{g_V\}_{V\in \V}$, where $g_V
=f_U|_V$ for $U=\kappa(V)$. The path connected components of open
subsets of $X$ are open (since $X$ is locally path connected) and
hence the family of path connected components of the sets $V\in
\V$ form an open cover of $X$. Using compactness, we may pass to
a finite subcover; therefore, without loss of generality we may
assume that $\V$ is finite and the sets $V\in \V$ are path
connected.

For any $V_1, V_2\in \V$ the function $ g_{V_1} - g_{V_2} :
V_1\cap V_2\to \R$ is locally constant. We claim that the set
$S_{V_1 V_2}\subset \R$ of real numbers $ g_{V_1}(x) - g_{V_2}(x)
\in \R$, where $x$ varies in $V_1\cap V_2$, is finite. Assume the
contrary, i.e. there exists an infinite sequence $x_n\in V_1\cap
V_2$, where $n=1, 2, \dots,$  such that
\begin{eqnarray}
g_{V_1}(x_n) - g_{V_2}(x_n)\not= g_{V_1}(x_m) - g_{V_2}(x_m)\quad\mbox{for}\quad n\not=m.\label{notequal}
\end{eqnarray}
By compactness we may assume that $x_n$ converges to a point
$x_\infty\in X$. Denote $U_1=\kappa(V_1), U_2=\kappa(V_2)$, where $U_1,
U_2\in \U$. Then $x_\infty$ belongs to $U_1\cap U_2$ and thus $x_n\in
U_1\cap U_2$ for all large $n$. Let $W\subset U_1\cap U_2$ denote
the path connected component of $x_\infty$ in $U_1\cap U_2$. Since $W$
is open and contains $x_\infty$, it follows that $x_n$ belongs to $W$
for all large enough $n$. The function $f_{U_1}(x) - f_{U_2}(x)$,
where $x\in U_1\cap U_2$, is continuous and locally constant hence
it is constant for $x\in W$. We obtain that
\begin{eqnarray*}
f_{U_1}(x_n) - f_{U_2}(x_n)= f_{U_1}(x_m) - f_{U_2}(x_m)
\end{eqnarray*}
for all large enough $n$ and $m$. But this contradicts (\ref{notequal}) since
$f_{U_i}|_{V_i}=g_{V_i}$ and therefore $f_{U_i}(x_n) =g_{V_i}(x_n)$ for $i=1, 2$ and all $n$.

The union
$$S=\bigcup_{V_1, V_2\in \V} S_{V_1 V_2}$$
of all subsets $S_{V_1V_2}$ is a
finite subset of the real line. We will show now that the subgroup
of $\R$ generated by $S$ contains the group of periods
$h_\xi(\pi_1(X,x_0))$. Let $\gamma:[0,1] \to X$ be an arbitrary
loop, $\gamma(0)=\gamma(1)=x_0$. We may find division points
$t_0=0< t_1<t_2<\dots < t_N=1$ and open sets $V_1, \dots V_N\in
\V$  such that $\gamma([t_{i-1},t_i]) \subset V_i$ for $i=1, 2,
\dots, N$. Then $\gamma(t_i)\in V_i\cap V_{i+1}$ for $i=1, 2,
\dots, N$, where we understand that $V_{N+1}=V_1$. According to
the definition of the line integral (see (\ref{eq1.2})) we have
\begin{eqnarray*}
h_\xi([\gamma]) \, =\,  \langle \xi, [\gamma]\rangle \, =\,
\int_\gamma\omega = \\
=\sum_{i=1}^N \left[g_{V_i}(\gamma(t_i)) -
g_{V_i}(\gamma(t_{i-1}))\right] =\\
=\sum_{i=1}^N \left[g_{V_i}(\gamma(t_i)) -
g_{V_{i+1}}(\gamma(t_{i}))\right],
\end{eqnarray*}
which shows that any period $h_\xi([\gamma])\in \R$ lies in the
subgroup generated by the finite set $S\subset \R$. This implies
that the group of periods is finitely generated and completes the
proof of the Lemma. \qed
\end{proof}

\begin{lemma}\label{lift} Assume that $X$ is connected and locally path connected.
Let $\omega$ be a continuous closed 1-form on $X$. Consider the covering map $p_\xi:\tilde X_\xi\to X$
determined by the \v Cech cohomology class $\xi=[\omega]\in \check H^1(X;\R)$ of $\omega$. Then $p_\xi^\ast(\omega) =dF,$
where $F:\tilde X_\xi\to \R$ is a continuous function.
\end{lemma}
\begin{proof} Note that the integral $\int_\gamma \tilde \omega =0$ vanishes for any closed loop $\gamma$ in $\tilde X_\xi$,
where $\tilde \omega $ denotes $p^\ast_\xi(\omega)$,
by the construction of the covering $\tilde X_\xi$. Define
$$F(\tilde x)=\int_{\tilde x_0}^{\tilde x}\tilde\omega,\qquad \tilde x \in \tilde X_\xi,$$
where $\tilde x_0\in \tilde X_\xi$
is a base point, and the integration is taken along an arbitrary path in $\tilde X_\xi$
connecting $\tilde x_0$ with $\tilde x$. It is easy to see that $F(\tilde x)$ is independent of the choice of the path, $F$
is continuous, and $dF=p^\ast_\xi(\omega)$. \qed
\end{proof}

The covering $\tilde X_\xi$ is now used to characterize the chain recurrent set $R_\xi$ as follows. 

\begin{prop}\label{projection}
Let $X$ be a connected and locally path connected compact metric
space. Given a continuous flow $\Phi: X\times \R\to X$ and a
cohomology class $\xi\in \check H^1(X;\R)$, consider the free
Abelian covering $p_\xi: \tilde X_\xi\to X$ associated with $\xi$
(see above) and the canonical lift $\tilde \Phi: \tilde
X_\xi\times \R \to \tilde X_\xi$ of the flow $\Phi$ to $\tilde X_\xi$.
Fix a metric $\tilde d$ on $\tilde X_\xi$, which is invariant
under the group of covering translations and such that the
projection $p_\xi$ is a local isometry. Then the chain recurrent
set $R_\xi=R_\xi(\Phi) \subset X$ coincides with $p_\xi(R(\tilde
\Phi))$, the image of the chain recurrent set $R(\tilde
\Phi)\subset \tilde X_\xi$ of the lifted flow under the
projection.
\end{prop}
\begin{proof}
Let $\varepsilon_0 >0$ be such that for any ball $\tilde B\subset
\tilde X_\xi$ of radius $\varepsilon_0$ (with respect to the metric
$\tilde d$) the following holds:

(1) $g\tilde B\cap \tilde B=\emptyset$ for any element $g\not=1$
of the group of covering translations of $\tilde X_\xi$, and

(2) the projection $p_\xi$ restricted to $\tilde B$ is an
isometry.

We may satisfy (1) since the group of covering transformations of
the covering $\tilde X_\xi$ acts properly discontinuously (see
\cite{Sp}, page 87). Note that $\xi|_B =0$, where $B=p_\xi(\tilde
B)$.

Let $\delta_0=\delta(\varepsilon_0)>0$ be the number given by Lemma \ref{lem4.1}.
Then the pair $(\varepsilon_0, \delta_0)$ is a scale of $\xi$ in the sense of Definition \ref{scale}.

Suppose that a point $\tilde x\in \tilde X_\xi$ belongs to the
chain recurrent set $R(\tilde \Phi)$. Then for any $\delta>0$ and
$T>0$ there exists a $(\delta, T)$-chain of the form $\tilde
x_0=\tilde x, \tilde x_1, \dots, \tilde x_{N-1}, \tilde x_N=\tilde
x$, $t_1, \dots, t_N\in \R$, such that $\tilde d(\tilde
x_{i-1}\cdot t_i, \tilde x_i)<\delta$ and $t_i\geq T$ for all $i=1, 2,
\dots, N$. We will assume below that $\delta < \delta_0$.
Projecting downstairs, we find a sequence
$$x_0=x= p_\xi(\tilde x),\,
x_1, \dots, x_N= p_\xi(\tilde x) =x\, \in \, X,$$ with
$x_i=p_\xi(\tilde x_i)$ satisfying $d(x_{i-1}\cdot t_i, x_i)<
\delta$ for $i=1, 2, \dots, N$. This sequence forms a $(\delta,
T)$-chain in $X$ which starts and ends at $x$. For any homology class $z\in H_1(X;\Z)$ associated with this chain
one has $\langle \xi, z\rangle =0$ since we can find a loop
representing this class admitting a lift to the covering
$\tilde X_\xi$. This shows that $p_\xi(R(\tilde \Phi))$ is
contained in $R_\xi$.

To prove the inverse inclusion, assume that $\tilde x \in \tilde
X$ is such that the point $x=p_\xi(\tilde x)\in X$ belongs to
$R_\xi$. Hence for any $\delta>0$ and $T>1$ we can find a
$(\delta, T)$-chain $x_0=x, x_1, \dots, x_N=x$, $t_i\in \R$, such
that $d(x_{i-1}\cdot t_i, x_i)<\delta$, $t_i\geq T$, and for any
associated homology class $z\in H_1(X;\Z)$ one has $\langle \xi,
z\rangle =0$. We will assume that $\delta <\delta_0$, where
$\delta_0$ is given as above. Choose continuous curves
$\sigma_i:[0,1]\to X$ such that $\sigma_i(0) =x_{i-1}\cdot t_i$
and $\sigma_i(1)=x_i$ for $i=1, 2, \dots, N$ and the image
$\sigma_i([0,1])$ is contained in a ball of radius $\epsilon_0$. 
The concatenation of the
parts of trajectories from $x_{i-1}$ to $x_{i-1}\cdot t_i$ and the
paths $\sigma_i$, where $i=1, 2, \dots, N$, forms a closed loop
$\gamma$, which starts and ends at $x$. This loop lifts to a closed loop in the
cover $\tilde X_\xi$ which starts and
ends at $\tilde x$ since the homology class $z=[\gamma]$ of the
loop satisfies $\langle \xi, z\rangle =0$. The lift $\tilde
\gamma$ of $\gamma$ is a concatenation of parts of trajectories of
the lifted flow $\tilde \Phi$ and the lifts $\tilde \sigma_i$ of
the paths $\sigma_i$, where $i=1, \dots, N$. We obtain points
$\tilde x_i\in \tilde X_\xi$, where $i=0, 1, \dots, N$, such that
$p_\xi(\tilde x_i)=x_i$ and $\tilde x_0=\tilde x=\tilde x_N$.
Besides, we have $\tilde \sigma_i(0)= \tilde x_{i-1}\cdot t_i$,
$\tilde \sigma_i(1) = \tilde x_i$ for $i=1, 2, \dots, N$. Since
each $\sigma_i$ lies in a ball of radius $\epsilon_0$ in $X$ it
follows from our assumption (1) above that each path $\tilde
\sigma_i$ lies in a ball $\tilde B\subset \tilde X_\xi$; from
assumption (2) we find that $\tilde d(\tilde x_{i-1}\cdot t_i,
\tilde x_i)<\delta$ for all $i=1, \dots, N$. Thus we have found a
$(\delta, T)$-chain in $\tilde X_\xi$ starting and ending at
$\tilde x$. This proves that ${p_\xi}^{-1}(R_\xi) \subset R(\tilde
\Phi)$, which is equivalent to $R_\xi\subset p_\xi(R(\tilde
\Phi))$. \qed
\end{proof}

\section{Proof of Theorem \ref{intromain}}

In this section we will prove our main Theorem \ref{intromain}.
The proof consists of two parts: the necessary conditions (easy)
and the sufficient conditions (more difficult).

\subsection{Necessary conditions}

If $\omega$ is a Lyapunov one-form for $(X, R_\xi)$ then 
by Definition \ref{introlyap}, 
$\omega|_{U}=df$ where $f$ is a continuous function 
defined on an open neighborhood $U\supset R_\xi$.
Hence the restriction of $\xi$ on $R_\xi$ vanishes, $\xi|_{R_\xi}=0$,
where $\xi\in \check H^1(X;\R)$ denotes the cohomology class of $\omega$. Thus condition (A) in Theorem
\ref{intromain} is necessary.
The following Proposition implies that condition (B) of Theorem \ref{intromain} is satisfied for any flow
admitting a Lyapunov one-form for $(\Phi, R_\xi)$ 
and having the property that the set $C_\xi=C_\xi(\Phi)$ is closed.

\begin{prop}
\label{prop4.5} Let $\Phi:X \x \bbr \to X$ be a continuous flow on
a compact, locally path connected, metric space $X$. Let $\omega$ be a Lyapunov one--form
for $(\Phi, Y)$, see Definition \ref{introlyap}, where
$Y\subset X$ is a closed flow-invariant subset.
Let $C\subset X$ be a closed, flow-invariant subset such that $Y\cap
C=\emptyset$. Then there exist numbers $\delta>0$ and $T>1$,
such that any homology class $z\in H_1(X;\Z)$
associated with any $(\delta, T)$-cycle $(x,t)$ with $x\in C$,
satisfies
$$\langle \xi, z\rangle \leq -1.$$
Here $\xi=[\omega] \in \check H^1(X;\R)$ denotes the cohomology
class of $\omega$.
\end{prop}
\begin{proof}
As $Y\cap C=\emptyset$ it follows that on $C$ the function $x\mapsto \int_x^{x\cdot 1}\omega\, <\, 0$
is continuous and negative.
Since $C$ is compact, there exists a positive
constant $c>0$, such that
\begin{eqnarray}\label{descend}
\int_x^{x\cdot 1}\omega <-c\quad\mbox{for all}\quad x\in C.
\end{eqnarray}

Let $(\varepsilon, \delta)$ be a scale of the cohomology class 
$\xi \in \check H^1(X;\R)$, see Definition \ref{scale}. 
Let $\eta>0$ be such that for
any continuous curve $\sigma: [0,1]\to X$, lying in a ball $B\subset X$ of
radius $\eps$ one has $|\int_\sigma \omega|<\eta$.
We define
\begin{eqnarray}\label{choice}
T=\left[{(1+\eta)}/{c}\right]+2,
\end{eqnarray}
where $[a]$ denotes the integer part of a number $a$. Unlike $\delta$, 
the number $T$ not only depends on the class $\xi$ but on the chosen representative $\omega$.
We will show that $\delta>0$ and $T>1$ satisfy our requirements.
Indeed, let $(x,t)$ be a
$(\delta, T)$-cycle with $x\in C$, and let $\gamma$ be a closed loop in $X$,
obtained by first following the trajectory $x\cdot\tau$, where
$\tau\in [0,t]$, and then returning from the end point $x\cdot t$ to $x$ along a
short path $\sigma$ lying in a  ball of radius $\eps$. Then we have
$$
\langle\xi, [\gamma] \rangle = \int_\gamma \omega = \int_x^{x\cdot t}\omega
+\int_\sigma\omega
$$
For the second integral we have $\int_\sigma\omega <\eta$ by
construction. Since $t \geq T$ and $T$ is an integer, we can estimate the first integral as
$$
\int_x^{x\cdot t}\omega = \sum_{i=1}^{T} \, \int_{x\cdot
(i-1)}^{x\cdot i}\omega +\int_{x\cdot T}^{x\cdot t}\omega < -Tc,
$$
where we have used $T$ times the inequality (\ref{descend}) and the estimate $\int_{x\cdot T}^{x\cdot t}\omega <0$.
By our choice of $T$, we have $Tc>1+\eta$ (cf. (\ref{choice})) and
so we see that
$
\langle \xi, z\rangle < -1
$
for any homology class $z\in H_1(X; \bbz)$ associated with any
$(\delta,T)$-cycle in $C$. \qed
\end{proof}

\subsection{Constructing a Lyapunov 1-form: the first step}

In the remainder of this section we prove the existence claim of Theorem \ref{intromain}.
The proof is split into several Lemmas. 
In a first step we construct a Lyapunov one-form for the flow restricted to $C_\xi$;
later on we will extend the obtained closed one-form to a Lyapunov one-form defined on the whole space $X$.

We start with the following combinatorial Lemma:
\begin{lemma}\label{lem5.2a}
Any nonempty word $w$ of arbitrary (finite) length
in an alphabet consisting of $L$ letters can be written as a
product (concatenation)
$$w=w_1w_2\dots w_l$$
of $l \leq L$ nonempty words, such that in each word $w_i$, the first and last letters coincide.
\end{lemma}

\begin{proof} We will use induction on $L$. For $L=1$ our claim is trivial.
We are left to prove the claim of the lemma assuming that it is
true for words in any alphabet consisting of less than $L$
letters. Consider a word $w$ in an alphabet with $L$ letters. Let
$a$ be the first letter of $w$. Finding the last appearance of $a$
in $w$, we may write $w=w_1 w'$, where $w_1$ starts and ends with
$a$ and $w'$ does not include $a$. We may now apply the induction
hypothesis to $w'$, which allows to write $w'=w_2 w_3\dots w_l$,
where $l\leq L$ and in each $w_i$ the first and the last letters
coincide. This clearly implies the lemma. \qed
\end{proof}

Let $\omega$ be an arbitrary, continuous closed one--form on $X$
representing a cohomology class $\xi\in \check H^1(X;\R)$. Our
final goal will be to modify $\omega$ so that at the end we
obtain a Lyapunov one--form for $(\Phi, R_\xi)$.

\begin{lemma}
\label{lem5.2} Under assumption (B) of Theorem
\ref{intromain} there exist $\mu>0$ and $\nu>0$ such that for
any $x\in C_\xi$ and $t\geq 0$ we have
\begin{eqnarray}
\int_x^{x\cdot t}\omega \leq -\mu t +\nu,\label{positive}
\end{eqnarray}
where the integral is calculated along the trajectory of the flow.
In particular,
$$\lim_{t\to +\infty}\int_x^{x\cdot t}\omega\, =\ -\infty$$
and the convergence is uniform with respect to
$x\in C_\xi$.
\end{lemma}
\begin{proof}
Let $\eps >0$ be such that for any ball $B\subset X$ of radius
$\varepsilon$ one has $\xi|_B=0$ (see Definition \ref{scale}) and for any continuous curve
$\sigma:[0,1]\to B$ holds $|\int_\sigma \omega|<1/2$. Let
$\delta>0$ be such that condition (B) of Theorem \ref{intromain}
holds for some $T>1$ and, additionally, for any points $x, y\in X$ with
$d(x,y)<\delta$ there is a continuous path $\sigma$ connecting $x$
and $y$ and lying in a ball of radius $\varepsilon$. By (B),
\begin{eqnarray}
\int_x^{x\cdot t}\omega \,  \, <\, \, -1/2, \label{half}
\end{eqnarray}
whenever $x\in C_\xi$, $t\geq T$ and $d(x,x\cdot t)<\delta$.

Using the compactness of $X$,
we find a constant $M>0$ such that for any point $x\in X$ and any time
$0\leq t\leq T$,
\begin{eqnarray}
\label{c}
\left|\int_x^{x\cdot t}\omega\right|<M.
\end{eqnarray}

Next we choose points $y_1, y_2, \dots, y_k$ in $X$, such that the
open balls of radius $\delta/2$ with centers at these points cover
$X$. Hence, for any $x\in X$ there exists an index $i\in \{1, \dots,
k\}$, such that $d(x, y_i)< \delta/2$.

Given $x\in C_\xi$ and $t\geq 0$, we consider the sequence of points
$$
x_j=x\cdot(jT) \, \in\,  C_\xi, \quad j=0,1, \dots, N,\quad\mbox{where}\quad N=[t/T].
$$
As explained above, for any $j=0, \dots, N$ there exists an index
$1\leq i_j\leq k$, such that
\begin{eqnarray}
d(x_j, y_{i_j})<\delta/2.\label{less}
\end{eqnarray}
Thus, any point $x\in C_\xi$ determines a sequence of indices
\begin{eqnarray}
i_0, i_1, \dots, i_N\in \{1, 2, \dots, k\},\label{index}
\end{eqnarray}
which, to a certain extent, encode the trajectory starting at $x$. If it happens that in
the sequence (\ref{index}) one has for some $r<s$, $i_r=i_s$, then the part
of the trajectory between $x_r=x\cdot (rT)$ and $x_s=x\cdot(sT)$
is a $(\delta, T)$-cycle (in view of (\ref{less})) and by (\ref{half})
\begin{equation}
\int_{x_r}^{x_s} \omega < -1/2.\label{janko}
\end{equation}
Let $1\leq m \leq k$ be an index which
appears in the sequence (\ref{index}) most often. Clearly, it must
appear at least $[(N+1)/k]$ times. Let $\alpha$ be the smallest
number with $i_\alpha =m$ and $\beta$ the largest number with
$i_\beta =m$, so that $0\leq \alpha <\beta \leq N$. Then, using
(\ref{janko}), we have
\begin{eqnarray}
\label{eq5.7}
\int_{x_{\alpha}}^{x_{\beta}}\omega \, \leq \,
-\, \frac{N+1-k}{2k}.
\end{eqnarray}
To complete the argument, we need to estimate the remaining
integrals $\int_x^{x_{\alpha}}\omega$ (corresponding to the
beginning of the trajectory) and $\int_{x_{\beta}}^{x\cdot
t}\omega$ (corresponding to the end of the trajectory).

View the sequence $i_0,i_1,\dots, i_{\alpha}$ as a word $w$ in the
alphabet $\{1, 2, \dots, k\}$ and apply Lemma \ref{lem5.2a}. As a
result we may split the sequence $w= i_0,i_1,\dots, i_{\alpha}$
into $l \leq k$ subsequences $w_1 , w_2, \dots, w_l$, each beginning
and ending with the same symbol. If $w_j=i_r, i_{r+1}, \dots, i_s$ is
one of the subsequences, where $r\leq s$, then $i_r=i_s$ and using
(\ref{janko}) we find $\int_{x_r}^{x_s}\omega \leq
0$. In other words, the integral corresponding to each subsequence
$w_j$ is nonpositive (in fact it is less than $-1/2$ if the
subsequence $w_j$ has more than one symbol).

Now we want to estimate the contribution of the integrals corresponding to the word breaks in
$w=w_1 w_2\dots w_l$. If $w_j$ ends with the symbol
$i_s$ and the following subsequence $w_{l+1}$ starts with
$i_{s+1}$ then we have $\int_{x_s}^{x_{s+1}}\omega \leq M$ (see
(\ref{c})). Thus, any word break contributes at most $M$ to the
integral.

The integral $\int_x^{x_\alpha}\omega$ is the sum of the
contributions corresponding to the words $w_j$ (which are all
nonpositive) and contributions of the word breaks (each is at most
$M$). Since there are $l-1 \leq k-1$ word breaks, we obtain
\begin{equation}
\label{eq5.8} \int_x^{x_{\alpha}}\omega \leq (l-1)M \leq (k-1)M.
\end{equation}

Similarly, $\int_{x_{\beta}}^{x_{N}}\omega
\leq (k-1)M. $ For the remaining integral we have $\int_{x_N}^{x\cdot t}\omega<M$, which again
follows from (\ref{c}), since $t-NT<T$.

Summing up, we finally obtain the estimate
\begin{eqnarray}\label{sum}
\int_x^{x\cdot t}\omega \, <\, (k-1)M - \frac{N+1-k}{2k} +(k-1)M
+M.
\end{eqnarray}
Hence (\ref{positive}) holds true with the constants 
$$\mu=\displaystyle{\frac{1}{2kT}} \quad \mbox{and}\quad
\nu=(2k-1)M+\frac{1}{2}.$$ 
\qed
\end{proof}

\begin{lemma} \label{closed}
Conditions {\rm (A)} and {\rm (B)} of Theorem \ref{intromain} imply that the set $C_\xi=R-R_\xi$ is closed.
\end{lemma}
\begin{proof} Since $\xi|_{R_\xi}=0$, we conclude that for any continuous closed one-form $\omega$ in the 
class $\xi$ the restriction $\omega|_{R_\xi}$ is the differential 
of a function and hence there exists a constant $C>0$ such that for any $x\in R_\xi$ and any $t>0$, 
\begin{eqnarray}
\left|\int_x^{x\cdot t}\omega\right|<C.\label{estimate1}
\end{eqnarray} 
Assume that the set $C_\xi$ is not closed, i.e. there exists a sequence of points $x_n\in C_\xi$ converging to a point 
$x_0\in R_\xi$. 
By Lemma \ref{lem5.2}, $\int\limits_{x_n}^{x_n\cdot t}\omega < -\mu t +\nu$. 
Taking $t=t_0= (\nu +2C)/\mu$, we obtain
\begin{eqnarray}\label{estimate}
\int_{x_n}^{x_n\cdot t_0}\omega \, <\, - 2C
\end{eqnarray}
for any $n=1, 2, \dots.$
Passing to the limit with respect to $n$ we find $\int\limits_{x_0}^{x_0\cdot t_0}\omega \, \leq\, - 2C$,
contradicting the estimate (\ref{estimate1}). \qed 
\end{proof}
\begin{lemma} \label{lem5.3}
Let $\omega$ be a continuous closed one-form on $X$ realizing a
class $\xi\in \check H^1(X;\R)$. Assume that conditions (A) and (B) of Theorem \ref{intromain} hold. 
Let $f: C_\xi\to \R$ be the function
defined by
\begin{eqnarray}
\label{deff}
f(x) := \sup_{t\geq 0}\, \int_x^{x\cdot t}\omega.
\end{eqnarray}
Then:
\begin{enumerate}[{\rm (i)}]
\item $f$ is well defined and continuous;
\item $\omega_1 =\omega|_{C_\xi} +df$ is a continuous closed one--form on
$C_\xi$ representing the cohomology class $\xi|_{C_\xi}$;
\item for any $x\in C_\xi$ and for any $t>0$,
\begin{equation}\label{eq5.10}
\int_x^{x\cdot t}\omega_1 \, \leq\,  0,
\end{equation}
i.e. $\omega_1$ is a Lyapunov one-form for the restricted flow
$\Phi|_{C_\xi}$ in a weak sense;
\item there exists a number $T>1$, such that for any $x\in C_\xi$ and any $t\geq T$,
\begin{equation}\label{eq5.11}
\int_x^{x\cdot t}\omega_1 \, \leq\,  -1.
\end{equation}
\end{enumerate}
\end{lemma}
\begin{proof} Let $t_0>0$ be the time such that $-\mu t_0+\nu =0$ with $\mu>0$ and $\nu>0$ as in Lemma \ref{lem5.2}.
 Then the supremum in (\ref{deff})
is achieved for $t\in [0,t_0]$, and hence we may write
$$
f(x) = \max_{0 \leq t \leq t_0} \,  \int_x^{x\cdot t}\omega.
$$
The continuity of $f$ now follows from the uniform
continuity of the integral with respect to $(x,t)\in C_\xi \x
[0,t_0]$.

Claim {\rm (ii)} follows from {\rm (i)}.
To prove {(iii)}, we find
\begin{eqnarray*}
\int_x^{x\cdot t} \omega_1 = \int_x^{x\cdot t} (\omega + df) &=&
\int_x^{x\cdot t} \omega +\left[f(x\cdot t) - f(x)\right] \\
&=& \int_x^{x\cdot t} \omega + \sup_{\tau\geq 0}\int_{x\cdot
t}^{x\cdot (t+\tau)}\!\!\! \omega \quad -\; \sup_{\tau\geq
0}\,\int_x^{x\cdot \tau}
\omega \\
&=& \sup_{\tau\geq t}\,\int_x^{x\cdot \tau} \omega \quad -\;
\sup_{\tau\geq 0}\,\int_x^{x\cdot \tau} \omega \, \leq\,  0.
\end{eqnarray*}

We next prove (iv). Let $M$ denote the maximal value of the continuous
function $f:C_\xi\to \R$ and $m$ its minimal value.
We apply Lemma \ref{lem5.2} and obtain 
\begin{eqnarray*}
\int_x^{x\cdot t} \omega_1 =
\int_x^{x\cdot t} \omega +\left[f(x\cdot t) - f(x)\right] \leq -\mu t +\nu +(M-m).
\end{eqnarray*}
Hence, choosing $T$ such that $-\mu T +\nu = -2 -(M-m)$, claim (iv) follows.
\qed
\end{proof}

\subsection{Second step: smoothing}\label{smoothing}

In this subsection we describe a procedure for smoothing
a continuous closed one-form along the flow, which will be used in
the proof of Theorem \ref{intromain} below. It is a modification of a
well-known method for continuous functions, see for example
\cite{Sc}. We use this construction to smooth the closed
one-form $\omega_1$ which is constructed in the proof of Lemma \ref{lem5.3}.

A function $f:X \to \bbr$ is said to be differentiable along
a continuous flow $\Phi:X \x \bbr \to X$ if the derivative $\frac{d}{dt} f(x
\cdot t)|_{t=0} $ exists for any $x\in X$. More generally, a
continuous closed one--form $\omega$ on $X$ is said to be {\it
differentiable with respect to the flow} $\Phi$ if the derivative
\begin{equation}
\dot\omega (x) :=
\frac{d}{dt} \left(\int_x^{x\cdot
t} \omega\right)|_{t=0}
\end{equation}
exists for any $x\in X$. In this case,
$\dot\omega :X\to \R$ is a function on $X$, which we call the
{\it derivative of $\omega$ with respect to the flow $\Phi$}.
If $\omega$ is represented as $\omega=\{\pHi_U\}_{U\in \U}$ with
respect to an open cover $\U$ of $X$ then for $x\in U$ and $t$
sufficiently small we have $ \int_x^{x\cdot t} \omega =
\pHi_{U}(x\cdot t) - \pHi_{U}(x)$ and we see that a continuous
closed one--form is differentiable with respect to the flow $\Phi$
if and only if the local defining functions $\pHi_{U}$ are.

\begin{lemma}
\label{lem1.3}
Assume that conditions (A), (B) of Theorem
\ref{intromain} hold. Then there exists a continuous closed one-form
$\omega_2$ on $X$ in class $\xi\in \check H^1(X;\R)$ with the following properties:
\begin{enumerate}[{\rm (i)}]
\item
$\omega_2$ is differentiable
with respect to the flow $\Phi$;
\item the derivative
$\dot\omega_2: X\to \R$ is a continuous function;
\item  for some $\sigma>0$ one has $\dot \omega_2(x) \leq -\sigma$ for all $x\in C_\xi$;
\item $\omega_2|_U =0$ for some open neighborhood $U\subset X$ of $R_\xi$.
\end{enumerate}
In particular, $\omega_2|_{C_\xi}$ is a
Lyapunov one-form for the flow $\Phi|_{C_\xi}$.
\end{lemma}

\begin{proof} Using assumptions (A), (B) and Lemma \ref{closed},
we find a closed neighborhood $V\subset X$ of $R_\xi$, such
that $\xi|_V=0$ and $V\cap C_\xi=\emptyset$. Here we use
the continuity property of the \v Cech cohomology theory, see \cite{CE},
Chapter 10, Theorem 3.1. 

Let $\omega_1$ be the closed one-form on $C_\xi$ given by Lemma
\ref{lem5.3}. Using the Tietze Extension Theorem (see Proposition
\ref{tietze}), we may find a closed one-form $\Omega_1$ on $X$,
such that $\Omega_1|_{C_\xi}=\omega_1$, $\Omega_1|_V=0$, and
$[\Omega_1]=\xi\in \check H^1(X;\R)$.

Consider the covering map
$p_\xi: \tilde X_\xi \to X$ corresponding to the \v Cech cohomology class $\xi$, see
subsection \ref{dynamics}. By Lemma \ref{lift} we have
$p^\ast_\xi(\Omega_1) =dF_1$, where $F_1: \tilde X_\xi\to \R$ is a
continuous function. Let $\tilde C_\xi$ denote the preimage
$p_\xi^{-1}(C_\xi)$. Then for any point $x\in \tilde C_\xi$ and $t\geq
0$, $F_1(x\cdot t) \leq F_1(x),$ by Lemma
\ref{lem5.3}, (iii). Moreover, statement (iv) of Lemma \ref{lem5.3}
implies that there exists $T>0$, such that
\begin{eqnarray}\label{less1}
F_1(x\cdot t) - F_1(x) \leq -1, \quad \mbox{for}\quad
x\in \tilde C_\xi, \, \, t\geq T.
\end{eqnarray}

Let $\varrho: \R \to [0,\infty)$ be a $C^\infty$-smooth function with the following properties:
\begin{enumerate}[{\rm (a)}]
\item the support of $\varrho$ is contained in the interval $[-T-1, T+1]$;
\item $\varrho|_{[-T,T]} ={\rm const} = \sigma >0$;
\item $\varrho(-t) =\varrho(t)$;
\item $\varrho'(t)\geq 0$ for $t\leq 0$;
\item $\int_\R \varrho(t)dt =1$.
\end{enumerate}
Using $\varrho$ we define $F_2: \tilde X_\xi \to \R$ by $F_2(x) = \int_\R F_1(x\cdot t) \varrho(t)dt$.
It is clear that $F_2$ is continuous. Since $F_2(x\cdot s) = \int_\R F_1(x\cdot t)\varrho(t-s)dt$, we see that
$F_2$ is differentiable with respect to the flow on $\tilde X_\xi$. If $x\in \tilde C_\xi$, we find, using (\ref{less1}) and the
properties of $\varrho$,
\begin{eqnarray}
\frac{dF_2(x\cdot s)}{ds}|_{s=0} &=& -\int_{-T-1}^{T+1} F_1(x\cdot t) \varrho'(t) dt \nonumber\\
&=&\int_{-T-1}^{-T} \left[ F_1(x\cdot (-t))- F_1(x\cdot t)\right]\cdot \varrho'(t)dt\label{derivative}\\
&\leq & -\int_{-T-1}^{-T}\varrho'(t)dt \, =\,  -\sigma.\nonumber
\end{eqnarray}
Let $G$ denote the group of covering transformations of the covering $\tilde X_\xi$.
Using the homomorphism of periods (\ref{periods}), one sees that 
the class $\xi$ determines a monomorphism $\alpha: G\to \R$, such that for any $x\in \tilde X_\xi$ and any $g\in G$
we have
\begin{eqnarray}\label{difference}
F_1(gx) - F_1(x) =\alpha(g).
\end{eqnarray}
Since $(gx)\cdot t=g(x\cdot t)$, we find
$F_1((gx)\cdot t) = F_1(g(x\cdot t)) = F_1(x\cdot t) + \alpha(g)$ and, multiplying by $\varrho(t)$ and integrating gives
\begin{eqnarray}\label{difference1}
F_2(gx) - F_2(x) =\alpha(g)
\end{eqnarray}
for any $x\in\tilde X_\xi$ and $g\in G$. Formula (\ref{difference1}) says that the action of the covering translations
changes $F_2$ by adding a constant, and therefore, $F_2$ determines a continuous closed one-form on $X$.
More precisely, $dF_2 = p_\xi^\ast(\omega_2)$ for some continuous closed one-form $\omega_2$ on $X$.
Since $F_2$ is differentiable with respect to the flow on $\tilde X_\xi$ and the derivative $\frac{d}{ds}F_2(x\cdot s)$ is continuous
(see (\ref{derivative})), the derivative
$\dot \omega_2: X\to \R$ is a well-defined, continuous function.
Clearly, as $\Omega_1$ vanishes on $V$, the form $\omega_2$ vanishes on the open set $U\subset X$ of points $x\in X$ with 
$x\cdot [-T-1,T+1]\subset V$. Since $R_\xi\subset U$, this proves (iv).  

Comparing (\ref{difference}) and (\ref{difference1}) we find that the function
$F_1-F_2: \tilde X_\xi \to \R$ is invariant under the covering translations.
Hence $F_1-F_2 = f\circ p_\xi$, where $f: X\to \R$ is a continuous function. Therefore $\Omega_1 -\omega_2 =df$, i.e.
$\omega_2$ lies in the cohomology class $\xi$. We know that $\omega_2$ is differentiable with respect to the flow and
$\dot \omega_2 \leq -\sigma<0$ on $C_\xi$.
\qed
\end{proof}

\subsection{Third step: extension}

Now we complete the proof of the existence claim of Theorem \ref{intromain}.

Let $L: X\to \R$ be a Lyapunov function for $(\Phi, R)$. Such a function
exists according to Theorem \ref{intro1} of C. Conley. We
apply the smoothing procedure from the previous subsection to $L$.
Namely, let $\rho: \R\to [0,\infty)$ be a $C^\infty$-smooth
function, such that $\supp (\rho) = [-1,1]$, $\int_\R \rho(t)dt
=1$, $\rho(-t) =\rho(t)$ and $\rho'(t) >0$ for all $t\in (-1, 0)$ and set
$$L_1(x) =\int_\R L(x\cdot t)\rho(t)dt.$$
We find (precisely as in the previous subsection) that the derivative $\dot
L_1(x)=\displaystyle{\frac{d}{ds}L_1(x\cdot s)|_{s=0}}$ exists and
is given by
\begin{eqnarray}
\dot L_1(x) = \int_{-1}^0\left[L(x\cdot (-t))-L(x\cdot t)\right]\rho'(t)dt.
\end{eqnarray}
This identity implies that $\dot L_1: X\to \R$ is a continuous
function and
$$\dot L_1(x)<0\quad\mbox{for any}\quad x\in X-R.$$

Let $\omega_2$ be the closed one-form on $X$
given by Lemma \ref{lem1.3}.
We will set
\begin{eqnarray}
\omega_3= \omega_2 +\lambda (dL_1),\label{final}
\end{eqnarray}
where $\lambda >0$ and $dL_1$ is the differential of the function
$L_1$ (see section \ref{prelim}). In view of the construction of $\omega_2$ and $L_1$, 
for any $\lambda $, the form $\omega_3$ is a continuous closed 1-form on $X$
representing the cohomology class $\xi$ and satisfies condition (L2) of
Definition
\ref{introlyap}.

We now show that for $\lambda$ large enough $\omega_3$ satisfies condition (L1) and hence it is
a Lyapunov one-form for $(\Phi, R_\xi)$.
Indeed, $\omega_3$ is differentiable along the flow and has the derivative 
$\dot \omega_3 = \dot\omega_2+\lambda\dot L_1$. 
By Lemma \ref{lem1.3}, $\dot \omega_2 <0$ on $C_\xi$. Hence
we may find an open neighborhood $W$ of $C_\xi$, so that $\dot
\omega_2<0$ on $W$. By claim (iv) of Lemma \ref{lem1.3}, $\dot \omega_2 =0$ vanishes on some open neighborhood
$U$ of $R_\xi$, whereas $\dot L_1<0$ on $U-R_\xi$. 
Hence we see that for any $\lambda >0$ the
inequality $\dot \omega_3 <0$ holds on $W$ and on $U-R_\xi$.
Finally we shall show that $\dot \omega_3 <0$ on $X-R_\xi$
for $\lambda >0$ sufficiently large.

The function
\begin{eqnarray}\label{function}
x\mapsto -\frac{\dot\omega_2(x)}{\dot L_1(x)}, \qquad x\in X-(U\cup W)
\end{eqnarray}
is well-defined and continuous (recall that $\dot L_1 < 0$ on
$X-R$). Since $X-(U\cup W)$ is compact, the function
(\ref{function}) is bounded. Choose $\lambda>0$ to be
larger than the maximum of (\ref{function}). Then $\dot
\omega_3(x) <0$ holds for all $x\in X-R_\xi$, as desired.

This completes the proof of Theorem \ref{intromain}. \qed

The arguments above prove the following, slightly stronger
statement:

\begin{corollary}\label{more}
Under the assumptions (A) and (B) of Theorem \ref{intromain},
there exists a continuous closed one-form $\omega$ on $X$ lying in the cohomology class
$[\omega]=\xi\in \check H^1(X;\R)$ which satisfies condition {\rm (L2)} and the following stronger version of condition
{\rm (L1)}: 
$\omega$ is differentiable with respect to
the flow $\Phi$ (in the sense explained in section \S
\ref{smoothing}), the derivative $\dot \omega: X\to \R$ is
continuous and $\dot \omega<0$ on $X-R_\xi$.
\end{corollary}
\qed

\section{Proof of Corollary \ref{introcor.2}}
\label{applic}

By Theorem \ref{intromain}, $\rm {(i)}$ implies $\rm {(ii)}$. 
Conversely, assume that there exists a Lyapunov one--form $\omega$ for
$(\Phi,\emptyset)$ representing the class $\xi$. Proposition \ref{prop4.5} shows that
condition $(B)$ of Theorem \ref{intromain} is satisfied. We are left to prove that $R_\xi=\emptyset$. 
Consider the covering $p_\xi: \tilde X_\xi \to X$ corresponding to the class $\xi$ (see \S \ref{dynamics}).
It is enough to show that the chain recurrent set $R(\tilde \Phi)$ of the lifted flow $\tilde \Phi$ in $\tilde X_\xi$ is empty
(see Proposition \ref{projection}). By Lemma \ref{lift}, $p_\xi^\ast(\omega) =dF$, where $F: \tilde X_\xi \to \R$
is a continuous function. By assumption (ii), $F(x\cdot t) < F(x)$ for all $x\in \tilde X_\xi$ and $t>0$. 
In particular, the function $\phi(x) = F(x\cdot 1)-F(x)$, defined on $x\in \tilde X_\xi$, 
is negative and invariant under the group of covering translations (see (\ref{difference})); 
hence it equals $\psi\circ p_\xi$, where $\psi: X\to \R$
is a continuous function on $X$. 
By the compactness of $X$ there exists $\sigma >0$ 
such that $\phi(x)<-\sigma$ for all $x\in \tilde X_\xi$. Choose a metric $d$ 
on $\tilde X_\xi$, which is invariant under the group
of covering translations. There exists $\delta>0$ such that for any $x,y\in \tilde X_\xi$ with $d(x,y)<\delta$ one has $|F(x)-F(y)|<\sigma/2$.
Now, assume that the points $x_0=x, x_1, \dots, x_N=x\in \tilde X_\xi$ and the numbers $t_1, \dots, t_N\in \R$ 
represent a $(\delta, T)$-chain with $T>1$ of the lifted flow $\tilde \Phi$ in $\tilde X_\xi$,
i.e. $t_i \geq T$ and $d(x_{i-1}\cdot t_i, x_i)<\delta$ for $i=1, \dots, N$. 
Then for any $i=1, 2, \dots, N$ we have
$$F(x_{i-1}\cdot t_i)-F(x_{i-1})< -\sigma\quad\mbox{and}\quad F(x_i) - F(x_{i-1}\cdot t_i) < \sigma /2,$$
which imply that $F(x_i)-F(x_{i-1})< -\sigma/2$ and hence $F(x_N)-F(x_0)<0$. 
The last inequality contradicts $x_0=x=x_N$. This proves that 
there are no closed $(\delta, T)$-chains in the covering $\tilde X_\xi$.

Thus we have shown that $\rm {(i)}$ and
$\rm {(ii)}$ are equivalent.

Now assume that $\rm {(ii)}$ holds and the class $\xi$ is integral, i.e. $\xi\in \check H^1(X;\Z)$.
As $X$ is
locally path connected and compact, it has finitely many
path connected components. Thus without loss of generality we may assume that
$X$ is path connected. Let $\omega$ be a Lyapunov one-form for $(\Phi, \emptyset)$ satisfying the properties (iii) and (iv) of Corollary \ref{more}.
Define a map $p:X \to S^1\subset \C$ by choosing a point $x_0\in X$ and setting
$$
p(x) = \exp\left[2\pi i\int_{x_0}^{x} \omega \right],
$$
where the line integral is taken along any path connecting $x_0$ with $x$.
Since $\xi$ is integral, the value $p(x)\in S^1$ does not depend on
the choice of the path. 
The function 
$$t\mapsto \arg(p(x\cdot t))= 2\pi  \int_{x_0}^{x\cdot t} \omega$$ 
is differentiable and the derivative 
$$
\displaystyle{\frac{d}{dt}\arg(p(x \cdot t))} <0
$$
is negative. 
The equality $\xi = p^*(\mu)$ is immediate from the
definition of $p$. 

It remains to prove that $p$ defines a locally trivial fibration.
Pick $\eta\in \R$ and let $K_{\eta}=p^{-1}(\exp(2\pi i\eta))\subset X$. For any $x\in K_\eta$, let $f_x: \R \to \R$ be a continuous
function such that $f_x(0)=\eta$ and $p(x\cdot t) = \exp(2\pi i f_x(t))$ for all $t\in \R$. The function $f_x$ is uniquely 
determined. It is 
differentiable and there exists $\varepsilon>0$ such that for all $x\in K_\eta$ and $t\in \R$, 
$$\frac{d}{dt}f_x(t)< -\varepsilon.$$
Let $g_x=f_x^{-1}$ be the inverse function. Define $G: K_\eta\times \R\to X$ by $G(x,t) = x\cdot g_x(t)$. 
Then $G$ is continuous and the diagram
\begin{eqnarray}
\begin{array}{ccc}
K_\eta \times \R & \stackrel G\lra & X\\ \\
\, \, e  \, \searrow &&\swarrow \, p\quad\\ \\
& S^1&
\end{array}
\end{eqnarray}
commutes, where $e(x,t) = \exp(2\pi i t)$. 
This proves that $p$ is a locally trivial fibration and that $K_\eta$ is a cross section of the flow $\Phi$.
\qed

\section{Proof of Proposition \ref{closing}}\label{proofclosing}

Suppose that the set $C_\xi\subset R=R(\Phi)$ is closed but condition (B) of Theorem \ref{intromain} is violated. 
Then there exists a sequence of
points $x_n\in C_\xi$ and numbers $t_n>0$ such that the distances $d(x_n, x_n\cdot t_n)$ tend to 0 as $t_n \to \infty$,
and 
\begin{eqnarray}
\langle \xi, z_n\rangle >-1,\label{new}
\end{eqnarray}
where $z_n\in H_1(X;\Z)$ denotes the homology class obtained by \lq\lq closing\rq\rq \, the trajectory
$x_n\cdot t$ for $t\in [0,t_n]$. 
Using the compactness of $C_\xi$ we may additionally assume that $x_n$ converges to a point $x\in C_\xi$. 
Since we assume that the class $\xi$ is integral, we may rewrite (\ref{new}) in the form 
$\langle \xi, z_n\rangle \geq 0$. 
Thus we obtain a 
closing sequence $(x_n,t_n)$ such that for any homology direction $\tilde z\in D_X$ associated with it, 
$\langle \xi, \tilde z\rangle \geq 0$, i.e. the condition of Fried \cite{Fr:1} also fails to hold. 

Conversely, we now show that the condition of Fried \cite{Fr:1}
holds assuming that condition (B) of Theorem \ref{intromain} is satisfied. 
Fix a norm $||\,\, ||$ on the vector space $H_1(X;\R)$. 
As $X$ is a polyhedron, there exists $\delta>0$ so that for any $\delta/2$-ball $B$ in $X$ the inclusion $B\to X$ is null-homotopic.
Further, there exists a constant $C>0$ such that for any homology class $z\in H_1(X;\Z)$ 
associated with a $(\delta, T)$-cycle $(x, t)$ in $X$ one has
\begin{eqnarray}
||z||\leq Ct.\label{again3}
\end{eqnarray} 
Let $\omega$ be a continuous 
closed 1-form in the class $\xi$. By Lemma \ref{lem5.2} there exist $\mu>0$ and $\nu>0$ such that 
\begin{eqnarray}
\int_x^{x\cdot t}\omega \leq -\mu t +\nu\label{again}
\end{eqnarray}
 for all $x\in C_\xi$ and $t>0$. Let $\eta>0$ be such that $\left|\int_\gamma \omega\right|<\eta$ for any curve lying in a ball 
of radius $\delta/2$. The estimate (\ref{again}) implies that 
\begin{eqnarray}\label{again1}
\langle \xi, z\rangle \leq -\mu t +\nu +\eta\end{eqnarray}
for any $(\delta, T)$ cycle $(x, t)$ with $x\in C_\xi$, where $z\in H_1(X;\Z)$ denotes the associated homology class. 
Since $\langle \xi, z\rangle \geq -c ||z||$, where $c>0$, we obtain that 
the homology class $z$ of any $(\delta, T)$-cycle $(x, t)$ with $x\in C_\xi$ satisfies
\begin{eqnarray}
||z|| \geq \frac{\mu}{c} \cdot t - \frac{\nu +\eta}{c}.\label{again2}
\end{eqnarray}
Now, let $(x_n, t_n)$ be a closing sequence (as defined in \S \ref{intro}), where $x_n\in C_\xi$, such that $x_n$ converges to a point $x\in C_\xi$ and  
$t_n\to \infty$. 
Let $z_n\in H_1(X;\Z)$ denote the homology class determined by closing $(x_n,t_n)$. 
Then (\ref{again2}) implies that $||z_n||\to \infty$. By (\ref{again3}), $-t \leq - \displaystyle{\frac{||z||}{C}}$, which when substituted 
into (\ref{again1}) leads to  
$$ \langle \xi, \frac{z_n}{||z_n||}\rangle \, \leq \, - \frac{\mu}{C}+\frac{\nu+\eta}{||z_n||}.$$ 
Therefore, we obtain 
for the homology direction $\displaystyle{\frac{z_n}{||z_n||}}\in D_X$ of the class $z_n$ the estimate
$$\langle \xi, \frac{z_n}{||z_n||}\rangle \, \leq \, - \frac{\mu}{2C}\, <\, 0,$$
if $n$ is large. 
This shows that Fried's condition \cite{Fr:1} is satisfied. \qed

\section{Examples}
\label{examples}

{\bf Example 1.} Here we describe  a class of examples of flows $\Phi: X\times \R\to X$,
for which there exists a cohomology class $\xi$ 
satisfying the conditions (A) and (B) of Theorem \ref{intromain}.

Let $M$ be a closed smooth manifold with a smooth vector field $v$. Let $\Psi: M\times \R\to M$ be the flow of $v$.
Assume that the chain recurrent set $R(\Psi)$ is a union of two disjoint closed sets $R(\Psi)=R_1\cup R_2$, 
where $R_1\cap R_2=\emptyset$. Out of this data we will construct a flow $\Phi$ on 
$$X=M\times S^1$$ 
such that 
$R_\xi(\Phi)=R_1\times S^0$, $C_\xi=R_2\times S^1$. Here $\xi\in H^1(X;\Z)$
denotes the cohomology class induced by the projection onto the circle $X\to S^1$ and $S^0\subset S^1$ is a two-point set.

Let $\theta\in [0,2\pi]$ denote the angle coordinate on the circle $S^1$. 
We will need two vector fields $w_1$ and $w_2$ on $S^1$, $w_1=cos (\theta)\cdot \frac{\partial}{\partial \theta}$ and
$w_2=\frac{\partial}{\partial \theta}$. The field $w_1$ has two zeros $\{p_1, p_2\}=S^0\subset S^1$ corresponding to
the angles $\theta=\pi/2$ and $\theta=3\pi/2$.

Let $f_i: M\to [0,1]$, where $i=1, 2$, be two smooth functions having disjoint supports and satisfying 
$f_1|_{R_1}=1$, $f_2|_{R_2}=1$.

Consider the flow $\Phi: X\times \R\to X$ determined by the vector field 
$$V\, = \, v+f_1w_1+f_2w_2.$$
Any trajectory of $V$ has the form $(\gamma(t), \theta(t))$, where $\dot\gamma(t)=v(\gamma(t))$, i.e.
$\gamma(t)$ is a trajectory of $v$. It follows that the chain recurrent set of $V$ is contained in $R(\Psi)\times S^1$. 
Over $R_1$ we have the vertical vector field $w_1$ along the circle 
which has two points $S^0\subset S^1$ as its chain recurrent set. 
Over $R_2$ we have the vertical vector field $w_2$ which has all of $S^1$ as the chain recurrent set. We see that 
$R_1\times S^0 = R_\xi(\Phi)$, $R_2\times S^1=C_\xi$. 
Hence $\xi|_{R_\xi}=0$ (and $C_\xi$ is closed). Clearly condition (B) of Theorem \ref{intromain} is satisfied as well.

{\bf Example 2.} 
Let $X=T^2$, thought of as $\bbr^2/\bbz^2$ with coordinates $x$ and $y$ on
$\bbr^2$. Any cohomology class $\xi\in H^1(T^2; \bbr)$ can be written as
$\xi = \mu [dx] + \nu [dy]$, where $dx$ and $dy$ are the standard coordinate 1-forms. 
We consider the flow of the following vector field 
$$V= f(x,y)\cdot \left(a \frac{\del}{\del x} + b \frac{\del}{\del y}\right),$$ 
where $b \neq 0$, $a/b\in \Q$ and $f:T^2 \to [0,1]$
is a smooth function vanishing at a single point $p\in T^2$. The chain recurrent set $R$ is the whole torus, $R=T^2$, while 
$$
R_\xi = \left\{ 
\begin{array}{cl}
T^2, & \text{ \rm if }\quad \mu a + \nu b = 0, \\ \\
f^{-1}(0)=\{p\}, & \text{ \rm otherwise}.
\end{array} \right.
$$
Assuming in addition that $\mu a + \nu b \neq 0$, the set $C_\xi = T^2-\{p\}$ is not closed. Nevertheless, a Lyapunov
one-form in the class $\xi\neq 0$ exists if and only if
$
\mu a + \nu b < 0.
$
In this case $\omega=\mu dx+\nu dy$ is such a Lyapunov 1-form.
This example shows that 
the existence of a Lyapunov one-form for $(\Phi, R_\xi)$
does not imply 
$C_\xi$ to be closed.

{\bf Example 3.}
Consider the standard irrational flow on the torus $X=T^2$, i.e. the flow of
the vector field $V= a \frac{\del}{\del x} + b \frac{\del}{\del y}$, where
$b \neq 0$ and $a/b\notin \Q$. Choose a cohomology class $\xi = \mu
[dx] + \nu [dy] \in H^1(X; \bbr)$ such that $\mu a + \nu b = 0$.
Then $R=X$ and $R_\xi=\emptyset$, but condition $(B)$ of Theorem \ref{intromain} is not satisfied
and so there is no Lyapunov one-form for $(\Phi, R_\xi)$ in the class $\xi$. 

This example shows that condition $(B)$ is not a consequence of the fact
that $C_\xi$ is closed.

\bibliographystyle{amsalpha}

\end{document}